\documentclass[a4paper,11pt]{article}

\pdfoutput=1
\usepackage[margin=1cm]{caption}

\usepackage{geometry}
\usepackage{authblk}
\usepackage{color}
\usepackage{times}
\usepackage{subfig}
\usepackage{bm}
\usepackage{graphicx}
\usepackage{graphics}
\DeclareGraphicsExtensions{.pdf,.png,.jpg}
\usepackage{amsbsy}
\usepackage{amsmath}
\usepackage{amsfonts}
\usepackage{amsthm}
\usepackage{float}
\usepackage{amssymb}
\usepackage{textcomp}
\usepackage{mathtools}
\usepackage[subnum]{cases}
\usepackage{bbold}
\usepackage{url}
\usepackage{framed}
\usepackage[colorlinks=true, citecolor=blue, urlcolor=blue]{hyperref}

\definecolor{purple(html/css)}{rgb}{0.5, 0.0, 0.5}

\newcommand{\cosec}{\operatorname{cosec}}
\newcommand{\me}{\mathrm{e}}
\newcommand{\mi}{\mathrm{i}}
\begin{document}

\title{A Few Finite Trigonometric Sums}
\author[1,2]{Chandan Datta \thanks{Email: \texttt{chandan@iopb.res.in}}} 
\author[1,2]{Pankaj Agrawal\thanks{Email: \texttt{agrawal@iopb.res.in}}}
\affil[1]{Institute of Physics, Sachivalaya Marg,
Bhubaneswar 751005, Odisha, India.}
\affil[2]{Homi Bhabha National Institute, Training School Complex, Anushakti Nagar, Mumbai 400085, India.}
\date{}

\maketitle

\begin{abstract}
    Finite trigonometric sums occur in various branches of physics,
    mathematics, and their applications. These sums may contain various powers of one or
    more trigonometric functions. Sums with one trigonometric function
    are known,  however sums with products of trigonometric functions
    can get complicated and may not have a simple expressions in a 
    number of cases.  Some of these sums have interesting
     properties and can have amazingly simple value. However,
    only some of them are available in literature. We obtain a number 
    of such sums using method of residues.

\end{abstract}

\section{Introduction}

      There is a venerable tradition of  computing finite sums of product
      of trigonometric functions in literature \cite{eisenstein,stern}. Such 
      sums occur while addressing many different problems in
      physics, or mathematics, or their applications.  Sums such as \cite{eisenstein}

\begin{eqnarray}
\sum_{j=1}^{d-1} \sin\bigg(\frac{2\pi mj}{d}\bigg)\cot\bigg(\frac{\pi j}{d}\bigg) = d - 2 m,
\end{eqnarray}
 where $d$ and $m$ are positive integers such that $m<d$, are known for a long time. However, if a small variation of arguments 
 of these functions is made, like arguments are affine functions, then
 these sums no longer remain easy to compute and are not available in various
 standard handbooks of mathematics \cite{gradshteyn}, including those specialize in series sum \cite{hansen, jolley, prudnikov}.
 As an example,  we may wish to compute
 \begin{eqnarray}
\sum_{j=1}^{d-1} \sin\bigg(\frac{2\pi mj}{d}+a\bigg) \cot\bigg(\frac{\pi j}{d}+\pi b\bigg). 
\end{eqnarray}
     There exist useful sum and difference formulas for sines and
     cosines that can be used, but such formulas do not exist for other
     trigonometric functions. In such cases,  there is a need to compute
     sums separately. For example, the above sum can be computed 
     for $ b = 0$, using existing results in the handbooks, but for
     $b\neq 0$, the sum is nontrivial. In this paper, our focus is on
     such sums.
  
    We have encountered such a sum in analyzing a Bell-type inequality
  for a system of two finite-dimensional subsystems.  We were computing
  Bell-Son-Lee-Kim (Bell-SLK) funcion for the general state of a bipartite quantum system.
  We encountered the
  following sum \cite{slk},
  \begin{equation}
 \sum\limits_{\alpha=0}^{d-1}\cos\Big(\frac{2\pi m}{d}(\alpha+\frac{1}{4})\Big)\cot\Big(\frac{\pi}{d}
(\alpha+\frac{1}{4})\Big),
\label{identity_3}
\end{equation}
where $m$ and $d$ are integers and satisfying $0<m<d$.
 We can use cosine sum rule, but there is no such rule for cotangent. So we cannot compute
 it using results given in standard handbooks.  This sum can be computed using a corollary, given below.
  Interestingly, this complicated looking sum
 has the value $d$, which is remarkably simple, and is independent of  $m$. Independence on $m$ is
 intriguing, and may have deeper mathematical meaning. Finite trigonometric sums have also appeared
 in the study of chiral Potts model \cite{chiral}, theory of Dirac operators \cite{fukumoto}, Dedekind sums \cite{fukuhara}, theory of determinants
 and permanents \cite{kittappa,minc}, and many other places.

      There are many techniques for computing these sums, e.g., use of generating functions,
      Fourier analysis, method of residues, etc. Sometimes, the same sum can be obtained by different
      methods, giving different looking results. Chu and Marini  \cite{chu} have used generating functions
      extensively to compute  the trigonometric sums. Berndt and Yeap \cite{berndt} have used method of residues.
      In this work, we will employ the method of residues \cite{berndt, brown}.
      We will start with a function with suitable singularity structure. The function and contour will
      be chosen in such a way that integration of the function over the contour
       gives the desired series and its sum.
      So  the main trick is to find suitable functions and computing residues at poles.
      
  In the next section, we have computed a few sums involving product of two different 
  trigonometric functions. In Section 3, we generalize the results to the sums of the
  product of more than two trigonometric functions. In Section 4, we conclude.
  
  \section{Products of Two Trigonometric Functions}

        In this section, we will compute the finite sums involving the product of sine with 
        powers of cotangent/cosecants and cosine with  powers of cotangent/cosecants.
        As a byproduct, we will also get a result involving tangent instead of cotangent.
        In computing these sums, we will use the  expansion
 
\begin{equation}\label{aexppoly}
\frac{1}{t\me^z-1}=\sum_{\nu=0}^{\infty}\frac{A_\nu(t)}{\nu!}z^\nu,
\end{equation}
where $A_\nu(t)$ is a function of $t$ and $t\neq 1$. The functions $A_{\nu}(t)$ can be written in terms of 
the so-called ``Apostol-Bernoulli numbers" $B_{\nu}(0,t)$\cite{apostol}. In fact
\begin{equation}
A_{\nu}(t)=\frac{B_{\nu+1}(0,t)}{\nu+1}.
\end{equation} 
The first few terms are $A_0(t)=\frac{1}{t-1}$, $A_1(t)=\frac{-t}{(t-1)^2}$, $A_2(t)=\frac{t+t^2}{(t-1)^3}$ and $A_3(t)=\frac{-(t+4t^2+t^3)}{(t-1)^4}$. 
We will also need to expand cotangent in a power series

\begin{equation}
\cot(\pi w)=\sum_{j=0}^{\infty}C_j\pi^{2j-1}w^{2j-1},
\end{equation}
where $w$ satisfies $0<|w|<\pi$ and 
\begin{equation}
C_j=\frac{(-1)^j2^{2j}B_{2j}}{(2j)!},
\end{equation}
where $B_{j}$ are the wellknown ``Bernoulli Numbers". 
 The first few $C_j$ are $C_0=1$, $C_1=-\frac{1}{3}$, $C_2=-\frac{1}{45}$ and $C_3=-\frac{2}{945}$. In our case, we need the expansion 
 \begin{equation}\label{cotpoly}
 \cot(\pi z+\pi b)=\sum_{j=0}^{\infty}C_j\pi^{2j-1}(z-1+b)^{2j-1},
\end{equation}     
 with the condition $0<|z-1+b|<\pi$.       
       Let us start with the following theorem. \\

\textbf{Theorem 1.1 :} If $m, n$ and $d$ denote positive integers with $m<d$ and $b\notin \mathbb{Z}/d$, then
\begin{eqnarray}\label{coscotn}
 e_n(d,m)&=& \sum_{j=0}^{d-1} \cos\bigg(\frac{2\pi mj}{d}\bigg)\cot^n\bigg(\frac{\pi j}{d}+\pi b\bigg)  \nonumber \\
    & = & -\sum \mi^{\mu+\nu+1}2^{\mu+\nu}\frac{m^\mu}{\mu !}\frac{d^{\nu+1}}{\nu !}\Big(t_1 A_\nu(t_2)-(-1)^{\mu+\nu}t^\prime_1 A_\nu(t^\prime_2)\Big)\nonumber\\&&\quad\quad D(j_1,j_2,\ldots,j_n).
 \end{eqnarray}
 Here the sum is over all nonnegative integers 
$j_1$,$\ldots$,$j_n$, $\mu$ and $\nu$ such that $2j_1+\cdots+2j_n+\mu+\nu=n-1$. We also have $t_1=\me^{-2\pi\mi mb}$, $t_2=\me^{-2\pi\mi db}$, $t^\prime_1=\me^{2\pi\mi mb}$ and $t^\prime_2=\me^{2\pi\mi db}$; here $b\notin \mathbb{Z}/d$ so that the trigonometric sum is well defined.

 Furthermore,
 \begin{equation}\label{dcotexp}
D(j_1,j_2,\ldots,j_{n})=\prod_{r=1}^{n}C_{j_r}.
\end{equation}
\\
\textbf{Proof :} We choose contour $C_R$ as a positively oriented indented rectangle with 
vertices at $\pm\mi R$ and $1\pm\mi R$. The contour has two semicircular indentations of 
radius $\epsilon \, (R>\epsilon)$ to the left of both $0$ and $1$ \cite{berndt}. Let us take the complex function as
\begin{equation}\label{coscotnfun}
f(z)=\frac{\me^{2\pi\mi mz }\cot^n(\pi z+\pi b)}{\me^{2\pi\mi dz}-1}
-\frac{\me^{-2\pi\mi mz }\cot^n(\pi z+\pi b)}{\me^{-2\pi\mi dz}-1}
\end{equation}
and consider $\frac{1}{2 \pi \mi}\int_C f(z)dz$.  Since $f(z)$ has period $1$, the integrals along the indented vertical sides of $C_R$ cancel. Since we have taken $m<d$, $f(z)$ tends to zero uniformly for $0\leqslant x\leqslant 1$ as $\vert y \vert\rightarrow\infty$. Hence, $\frac{1}{2 \pi \mi}\int_{C} f(z)dz=0$.  We can now calculate the contour integral using
Cauchy's residue theorem. The function $f(z)$ has poles at a number of points. To start with,
$f(z)$ has a simple pole at $z=0$, with residue 
\begin{equation}\label{coscotres0}
\mbox{Res}(f,0)=\frac{1}{\pi \mi d} \cot^n(b\pi).
\end{equation}
The function $f(z)$ also has simple poles at $z=\frac{j}{d}$, with $1\leqslant j\leqslant d-1$. The corresponding residues at these points are 
\begin{equation}\label{coscotnresj}
\mbox{Res}\bigg(f,\frac{j}{d}\bigg)=\frac{1}{\pi\mi d}\cos\bigg(\frac{2\pi mj}{d}\bigg)\cot^n\bigg(\frac{\pi j}{d}+\pi b\bigg).
\end{equation} 
In addition the function $f(z)$ has a pole of order $n$ at $z=-b+1$. Using equations (\ref{aexppoly}) and (\ref{cotpoly}),
 we can write
\begin{eqnarray}\label{coscotnfunexp}
f(z)&=&t_1\sum_{\mu=0}^{\infty}\frac{(2\pi \mi m)^\mu}{\mu!}(z+b-1)^\mu\Bigg(\sum_{j=0}^{\infty}C_j\pi^{2j-1}(z-1+b)^{2j-1}\Bigg)^{n}\nonumber\\&&\quad\sum_{\nu=0}^{\infty}
\frac{(2\pi \mi d)^\nu}{\nu!}A_\nu(t_2) (z+b-1)^\nu\nonumber\\&& -t^\prime_1\sum_{\mu=0}^{\infty}\frac{(-2\pi \mi m)^\mu}{\mu!}(z+b-1)^\mu\Bigg(\sum_{j=0}^{\infty}C_j\pi^{2j-1}(z-1+b)^{2j-1}\Bigg)^{n}\nonumber\\&&\quad\sum_{\nu=0}^{\infty}
\frac{(-2\pi \mi d)^\nu}{\nu!}A_\nu(t^\prime_2)(z+b-1)^\nu.\
\end{eqnarray}
Then after few steps of straightforward calculation, one can show that,
\begin{equation}\label{coscotnresb}
\mbox{Res}(f,-b+1)=\sum \mi^{\mu+\nu}2^{\mu+\nu}\frac{m^\mu}{\mu !}\frac{d^\nu}{\nu !}\big(\frac{t_1}{\pi} A_\nu(t_2)- (-1)^{\mu+\nu}\frac{t^\prime_1}{\pi} A_\nu(t^\prime_2)\big)D(j_1,j_2,\ldots,j_n).
\end{equation}
Here the sum is over all nonnegative integers 
 $j_1$,$\ldots$,$j_n$, $\mu$ and $\nu$ such that $2j_1+\cdots+2j_n+\mu+\nu=n-1$.
Using (\ref{coscotres0}), (\ref{coscotnresj}), (\ref{coscotnresb}) and applying residue theorem we can  obtain the sum (\ref{coscotn}).
\\
\\
\textbf{Corollary 1.2 :} Let $m$ and $d$ be positive integers such that $m<d$ and $b\notin \mathbb{Z}/d$. Then 
\begin{equation}\label{coscotsum}
e_1(d,m)=d \cos[(2m-d)b\pi]\cosec (bd\pi).
\end{equation}
\\
\textbf{Proof :} Put $n=1$ in Theorem 1.1. Using the values $A_0(t)$ and $C_0$, one can easily see this.
\\
\\
 We will now consider sums involving sine and cotangents. We will need to modify our functions suitably. \\
    
\textbf{Theorem 1.3 :} Let $m$,$n$ and $d$ denote positive integers with $m<d$ and $b\notin \mathbb{Z}/d$. Then
\begin{eqnarray}\label{sincotn}
 g_n(d,m)&=&-\sum \mi^{\mu+\nu}2^{\mu+\nu}\frac{m^\mu}{\mu !}\frac{d^{\nu+1}}{\nu !}\Big(t_1 A_\nu(t_2)+(-1)^{\mu+\nu}t^\prime_1 A_\nu(t^\prime_2)\Big)\nonumber\\&&\quad\quad D(j_1,j_2,\ldots,j_n).
 \end{eqnarray}
 We have already defined all the terms and conditions in Theorem 1.1. Here
\begin{equation}\label{sincotnsum}
g_n(d,m)=\sum_{j=1}^{d-1} \sin\bigg(\frac{2\pi mj}{d}\bigg)\cot^n\bigg(\frac{\pi j}{d}+\pi b\bigg).
\end{equation}
\\
\textbf{Proof :} Our contour will be the same as in Theorem 1.1. Let us take the complex function as
\begin{equation}\label{sincotnfun}
f(z)=\frac{\me^{2\pi\mi mz }\cot^n(\pi z+\pi b)}{\me^{2\pi\mi dz}-1}
+\frac{\me^{-2\pi\mi mz }\cot^n(\pi z+\pi b)}{\me^{-2\pi\mi dz}-1}
\end{equation}
and consider $\frac{1}{2 \pi \mi}\int_C f(z)dz$. As before, $\int_{C} f(z)dz=0$. The pole structure of this function 
is same as in Theorem 1.1. However, this time residue is zero at $z=0$, so we take the sum from $j=1$. Since $\sin x$ vanishes at $x=0$, we can  take the sum from $j=0$.
The function $f(z)$ has simple poles at $z=\frac{j}{d}$, with $1\leqslant j\leqslant d-1$. The corresponding residues at these points are 
\begin{equation}\label{sincotnresj}
\mbox{Res}\bigg(f,\frac{j}{d}\bigg)=\frac{1}{\pi d}\sin\bigg(\frac{2\pi mj}{d}\bigg)\cot^n\bigg(\frac{\pi j}{d}+\pi b\bigg).
\end{equation} 
The function $f(z)$ also has a pole of order $n$ at $z=-b+1$. Using (\ref{aexppoly}) and (\ref{cotpoly}), as in the case of
last theorem,  we can obtain after a  few steps of straight forward calculation,

\begin{equation}\label{sincotnresb}
\mbox{Res}(f,-b+1)=\sum \mi^{\mu+\nu}2^{\mu+\nu}\frac{m^\mu}{\mu !}\frac{d^\nu}{\nu !}\big(\frac{t_1}{\pi} A_\nu(t_2)+ (-1)^{\mu+\nu}\frac{t^\prime_1}{\pi} A_\nu(t^\prime_2)\big)D(j_1,j_2,\ldots,j_n).
\end{equation}
Using (\ref{sincotnresj}), (\ref{sincotnresb}) and applying residue theorem we can easily obtain (\ref{sincotn}).
\\
\\
\textbf{Corollary 1.4 :} Let $m$ and $d$ be positive integers such that $m<d$ and $b\notin \mathbb{Z}/d$. Then 
\begin{equation}\label{sincotsum}
g_1(d,m)=-d \sin[(2m-d)b\pi]\cosec (bd\pi).
\end{equation}
\\
\textbf{Proof :} Put $n=1$ in Theorem 1.3. Using the values $A_0(t)$ and $C_0$, we can obtain this.
\\
\\
 We will now consider sums involving sine and cosecants. We will do for even and odd powers of  cosecants
     separately. \\

\textbf{Theorem 1.5 :} Let $m$,$n$ and $d$ denote positive integers with $m<d$ and $b \notin \mathbb{Z}/d$. Then
\begin{eqnarray}\label{sincosec2n}
 h_n(d,m)&=&-\sum \mi^{\mu+\nu}2^{\mu+\nu}\frac{m^\mu}{\mu !}\frac{d^{\nu+1}}{\nu !}\Big(t_1 A_\nu(t_2)+(-1)^{\mu+\nu}t^\prime_1 A_\nu(t^\prime_2)\Big)\nonumber\\&&\quad\quad F(j_1,j_2,\ldots,j_{2n}),
 \end{eqnarray}
where the sum is over all nonnegative integers $j_1$,$\ldots$,$j_{2n}$, $\mu$ and $\nu$ such that $2j_1+\cdots+2j_{2n}+\mu+\nu=2n-1$ and $\mu+\nu$ must be odd. 
All other terms have already been defined in Theorem 1.1 except for $F(j_1,j_2,\ldots,j_{2n})$. Here
\begin{equation}\label{sincosec2nsum}
h_n(d,m)=\sum_{j=1}^{d-1} \sin\bigg(\frac{2\pi mj}{d}\bigg)\cosec^{2n}\bigg(\frac{\pi j}{d}+\pi b\bigg).
\end{equation}

The function $F(j_1,j_2,\ldots,j_{2n})$ is defined using the expansion
\begin{equation}\label{cosecpoly}
\cosec (\pi z+\pi b)=-\cosec(\pi z+\pi b-\pi)=\sum_{j=0}^{\infty}G_j\pi^{2j-1}(z-1+b)^{2j-1}, 
\end{equation}
with $0<|z-1+b|<\pi$.
Here 
\begin{equation}
G_j=\frac{(-1)^j2(2^{2j-1}-1)B_{2j}}{(2j)!},
\end{equation}
where $B_j$ is Bernoulli number. The first few $G_j$ are $G_0=-1$, $G_1=-\frac{1}{6}$, $G_2=-\frac{7}{360}$ and $G_3=-\frac{31}{15120}$. Then
\begin{equation}\label{fbcosecexp}
F(j_1,j_2,\ldots,j_{2n})=\prod_{r=1}^{2n}G_{j_r}.
\end{equation}
\\
\textbf{Proof :} Let's take the complex function as
\begin{equation}\label{sincosec2nfun}
f(z)=\frac{\me^{2\pi\mi mz }\cosec^{2n}(\pi z+\pi b)}{\me^{2\pi\mi dz}-1}
+\frac{\me^{-2\pi\mi mz }\cosec^{2n}(\pi z+\pi b)}{\me^{-2\pi\mi dz}-1}.
\end{equation}
We can use  the same contour and follow the same procedure as in Theorem 1.1 and Theorem 1.3.
The function $f(z)$ has simple poles at $z=\frac{j}{d}$, with $1\leqslant j\leqslant d-1$  with residues \\
\begin{equation}\label{sincosec2nresj}
\mbox{Res}\bigg(f,\frac{j}{d}\bigg)=\frac{1}{\pi d}\sin\bigg(\frac{2\pi mj}{d}\bigg)\cosec^{2n}\bigg(\frac{\pi j}{d}+\pi b\bigg).
\end{equation}
The function $f(z)$ also has a pole of order $2n$ at $z=-b+1$. Using (\ref{aexppoly}) and (\ref{cosecpoly}) we can write
\begin{eqnarray}\label{sincosec2nfunexp}
f(z)&=&t_1\sum_{\mu=0}^{\infty}\frac{(2\pi \mi m)^\mu}{\mu!}(z+b-1)^\mu\Bigg(\sum_{j=0}^{\infty}B_j\pi^{2j-1}(z-1+b)^{2j-1}\Bigg)^{2n}\nonumber\\&&\quad\sum_{\nu=0}^{\infty}
\frac{(2\pi \mi d)^\nu}{\nu!}A_\nu(t_2)(z+b-1)^\nu\nonumber\\&& +t^\prime_1\sum_{\mu=0}^{\infty}\frac{(-2\pi \mi m)^\mu}{\mu!}(z+b-1)^\mu\Bigg(\sum_{j=0}^{\infty}B_j\pi^{2j-1}(z-1+b)^{2j-1}\Bigg)^{2n}\nonumber\\&&\quad\sum_{\nu=0}^{\infty}
\frac{(-2\pi \mi d)^\nu}{\nu!}A_\nu(t^\prime_2)(z+b-1)^\nu.\
\end{eqnarray}
After a few steps of straightforward calculation, one obtains
\begin{equation}\label{sincosec2nresb}
\mbox{Res}(f,-b+1)=\sum \mi^{\mu+\nu}2^{\mu+\nu}\frac{m^\mu}{\mu !}\frac{d^\nu}{\nu !}\big(\frac{t_1}{\pi} A_\nu(t_2)+ (-1)^{\mu + \nu} \frac{t^\prime_1}{\pi} A_\nu(t^\prime_2)\big)F(j_1,j_2,\ldots,j_n).
\end{equation}
Using (\ref{sincosec2nresj}), (\ref{sincosec2nresb}) and applying residue theorem we can easily obtain (\ref{sincosec2n}).
\\
\textbf{Corollary 1.6 :} Let $m$ and $d$ be positive integers such that $m<d$ and $b \notin \mathbb{Z}/d$. Then 
\begin{equation}\label{sincosec2sum}
h_1(d,m)=d \cosec^2\big((b-1)d\pi\big)\Big[m \sin\big(2(b-1)(d-m)\pi\big) -(d-m) \sin\big(2(b-1)m\pi\big)\Big].
\end{equation}
\\
\textbf{Proof :} Put $n=1$ in Theorem 1.5. Using the values $A_0(t)$, $A_1(t)$ and $B_0$, one can show this.
\\

We will next  consider a sum involving cosine and even powers of the cosecant. We will have
to only change the function slightly. \\

\textbf{Theorem 1.7 :} Let $m$,$n$ and $d$ denote positive integers with $m<d$ and $b\notin \mathbb{Z}/d$. Then
\begin{eqnarray}\label{coscosec2n}
 k_n(d,m)&=&-\sum \mi^{\mu+\nu+1}2^{\mu+\nu}\frac{m^\mu}{\mu !}\frac{d^{\nu+1}}{\nu !}\Big(t_1 A_\nu(t_2)-(-1)^{\mu+\nu}t^\prime_1 A_\nu(t^\prime_2)\Big)\nonumber\\&&\quad\quad F(j_1,j_2,\ldots,j_{2n}),
 \end{eqnarray}
where the sum is over all nonnegative integers $j_1$,$\ldots$,$j_{2n}$, $\mu$ and $\nu$ such that $2j_1+\cdots+2j_{2n}+\mu+\nu=2n-1$ and $\mu+\nu+1$ must be even. Here
\begin{equation}\label{coscosec2nsum}
k_n(d,m)=\sum_{j=0}^{d-1} \cos\bigg(\frac{2\pi mj}{d}\bigg)\cosec^{2n}\bigg(\frac{\pi j}{d}+\pi b\bigg).
\end{equation}
\\
\textbf{Proof :} Let's take the complex function as
\begin{equation}\label{coscosec2nfun}
f(z)=\frac{\me^{2\pi\mi mz }\cosec^{2n}(\pi z+\pi b)}{\me^{2\pi\mi dz}-1}
-\frac{\me^{-2\pi\mi mz }\cosec^{2n}(\pi z+\pi b)}{\me^{-2\pi\mi dz}-1}.
\end{equation}
Then applying the same procedure as in Theorem 1.5, we can obtain (\ref{coscosec2n}).
\\
\\
\textbf{Corollary 1.8 :} Let $m$ and $d$ be positive integers such that $m<d$ and $b\notin \mathbb{Z}/d$. Then 
\begin{equation}\label{coscosec2sum}
k_1(d,m)=d \cosec^2\big((b-1)d\pi\big)\Big[m \cos\big(2(b-1)(d-m)\pi\big) +(d-m) \cos\big(2(b-1)m\pi\big)\Big].
\end{equation}
\\
\textbf{Proof :} Put $n=1$ in Theorem 1.7. Using the values $A_0(t)$, $A_1(t)$ and $B_0$, this result can be obtained.
\\

   Let us now consider odd powers of cosecants. In this case to keep the periodicity of the function $f(z)$, we shall
   have to use the argument of sine as $\frac{\pi mj}{d}$, instead of $\frac{2 \pi mj}{d}$ as earlier. \\

\textbf{Theorem 1.9 :} If $m$ is an odd positive integer and $n$, $d$ are positive integers with $m<d$, then
\begin{eqnarray}\label{sincosec2nodd}
 l_n(d,m)&=&-\sum \mi^{\mu+\nu}2^{\nu}\frac{m^\mu}{\mu !}\frac{d^{\nu+1}}{\nu !}\Big(p_1 A_\nu(p_2)+(-1)^{\mu+\nu}p^\prime_1 A_\nu(p^\prime_2)\Big)\nonumber\\&&\quad\quad F(j_1,j_2,\ldots,j_{2n-1}),
 \end{eqnarray}
where the sum is over all nonnegative integers $j_1$,$\ldots$,$j_{2n-1}$, $\mu$ and $\nu$ such that $2j_1+\cdots+2j_{2n-1}+\mu+\nu=2n-2$ and $\mu+\nu$ must be even. Here $p_1=\me^{-\pi\mi m(b-1)}$, $p_2=\me^{-2\pi\mi db}$, $p^\prime_1=\me^{\pi\mi m(b-1)}$ and $p^\prime_2=\me^{2\pi\mi db}$; $b\notin \mathbb{Z}/d$ such that the trigonometric sum is well defined.
\begin{equation}\label{sincosec2noddsum}
l_n(d,m)=\sum_{j=1}^{d-1} \sin\bigg(\frac{\pi mj}{d}\bigg)\cosec^{2n-1}\bigg(\frac{\pi j}{d}+\pi b\bigg).
\end{equation} 
\\
\textbf{Proof :} Let's take the complex function as
\begin{equation}\label{sincosec2noddfun}
f(z)=\frac{\me^{\pi\mi mz }\cosec^{2n-1}(\pi z+\pi b)}{\me^{2\pi\mi dz}-1}
+\frac{\me^{-\pi\mi mz }\cosec^{2n-1}(\pi z+\pi b)}{\me^{-2\pi\mi dz}-1}.
\end{equation}
Again proceeding in the same way as before we can obtain (\ref{sincosec2nodd}).
\\
\\
\textbf{Corollary 1.10 :} If $m$ is an odd positive integer and $n$, $d$ are positive integers with $m<d$ and $b\notin \mathbb{Z}/d$, then
\begin{equation}\label{sincosecsum}
l_1(d,m)=-d \sin\big((m-d)b\pi\big)\cosec(\pi db).
\end{equation}
\\
\textbf{Proof :} By putting $n=1$ in Theorem 1.9 and taking the values of  $A_0(p)$, and $B_0$, one can easily see this.
\\
   Again as before, in the above sum sine can be replaced with cosine, by modifying $f(z)$ slightly.  \\

\textbf{Theorem 1.11 :} If $m$ is an positive odd integer and $n$, $d$ are positive integers with $m<d$ and $b\notin\mathbb{Z}/d$, then
\begin{eqnarray}\label{coscosec2nodd}
 q_n(d,m)&=&-\sum \mi^{\mu+\nu+1}2^{\nu}\frac{m^\mu}{\mu !}\frac{d^{\nu+1}}{\nu !}\Big(p_1 A_\nu(p_2)-(-1)^{\mu+\nu}p^\prime_1 A_\nu(p^\prime_2)\Big)\nonumber\\&&\quad\quad F(j_1,j_2,\ldots,j_{2n-1}),
 \end{eqnarray}
where the sum is over all nonnegative integers $j_1$,$\ldots$,$j_{2n-1}$, $\mu$ and $\nu$ such that $2j_1+\cdots+2j_{2n-1}+\mu+\nu=2n-2$ and 
$\mu+\nu$ must be even. Here,
\begin{equation}\label{coscosec2noddsum}
q_n(d,m)=\sum_{j=0}^{d-1} \cos\bigg(\frac{\pi mj}{d}\bigg)\cosec^{2n-1}\bigg(\frac{\pi j}{d}+\pi b\bigg).
\end{equation} 
\\
\textbf{Proof :} Let's take the complex function as
\begin{equation}\label{coscosec2noddfun}
f(z)=\frac{\me^{\pi\mi mz }\cosec^{2n-1}(\pi z+\pi b)}{\me^{2\pi\mi dz}-1}
-\frac{\me^{-\pi\mi mz }\cosec^{2n-1}(\pi z+\pi b)}{\me^{-2\pi\mi dz}-1}.
\end{equation}
Again proceeding in the same way as before we can obtain (\ref{coscosec2nodd}).
\\
\\
\textbf{Corollary 1.12 :} If $m$ is an odd positive integer and $n$, $d$ are positive integers with $m<d$ and $b\notin\mathbb{Z}/d$, then
\begin{equation}\label{coscosecsum}
q_1(d,m)=d \cos\big((m-d)b\pi\big)\cosec(\pi db).
\end{equation}
\\
\textbf{Proof :} By putting $n=1$ in Theorem 1.11 and taking the values of  $A_0(p)$, and $B_0$, we can easily obtain this.
\\

If we take tangent instead of cotangent in the sum of corollary 1.2, then we will get 
\begin{equation}\label{costansum}
  \sum_{j=0}^{d-1} \cos\bigg(\frac{2\pi mj}{d}\bigg)\tan\bigg(\frac{\pi j}{d}+\pi b\bigg)  =
    \begin{cases}
      (-1)^{m+1}d \cos[(2m-d)b\pi]\cosec (bd\pi); & \text{if}\ d=\text{even}, \\
     (-1)^{m-d}d \sin[(2m-d)b\pi]\sec (bd\pi) ; & \text{if}\ d=\text{odd}.
    \end{cases}
  \end{equation}

Similarly, if in the sum of corollary 1.4, instead of cotangent, we take tangent then
\begin{equation}\label{sintansum}
  \sum_{j=1}^{d-1} \sin\bigg(\frac{2\pi mj}{d}\bigg)\tan\bigg(\frac{\pi j}{d}+\pi b\bigg)  =
    \begin{cases}(-1)^{m}d \sin[(2m-d)b\pi]\cosec (bd\pi); & \text{if}\ d=\text{even}, \\
   (-1)^{m-d}d \cos[(2m-d)b\pi]\sec (bd\pi) ; & \text{if}\ d=\text{odd}.
    \end{cases}
  \end{equation}

For even $d$  the sums in (\ref{costansum}) and (\ref{sintansum}) are not defined when $b \in\mathbb{Z}/d$, and for odd $d$ the sums are not defined when $b \in \mathbb{O}/2d$, where $\mathbb{O}$ is the set of odd integers.

Many more sums of the product of two trigonometric functions can be obtained with the same contour, with a change of function.
As an example, if we change slightly the arguments of the exponentials of the complex functions $f(z)$, we can get different kind of sums. 
To illustrate this, we can take the factor of $4$ instead of $2$ in the exponentials of the function $f(z)$ in Theorems 1.1 and 1.3.
Following the same procedure as used there, we get
\begin{equation}\label{coscotsum2d}
\sum_{j=0}^{2d-1} \cos\bigg(\frac{2\pi mj}{d}\bigg)\cot\bigg(\frac{\pi j}{2d}+\pi b\bigg)=2d \cos[(2m-d)2b\pi]\cosec (2bd\pi)
\end{equation}
and 
\begin{equation}\label{sincotsum2d}
\sum_{j=1}^{2d-1} \sin\bigg(\frac{2\pi mj}{d}\bigg)\cot\bigg(\frac{\pi j}{2d}+\pi b\bigg)=-2d \sin[(2m-d)2b\pi]\cosec (2bd\pi)
\end{equation}
respectively and $b\notin \mathbb{Z}/2d$ such that the sum is well defined.


      These are just a sample of sums that one can obtain with the contour described in Theorem 1.1, but choosing a variety
      of integrands.
      
  \section{Products of More than Two Trigonometric Functions}

We shall now consider a few sums involving the products of three trigonometric functions.
We shall consider only those functions for which integral along
 the two sides of the contour  cancel each other. Let us consider the following function
\begin{equation}\label{coscoseccosfun}
f(z)=\bigg[\frac{\me^{2\pi\mi mz }}{\me^{2\pi\mi dz}-1}-\frac{\me^{-2\pi\mi mz }}{\me^{-2\pi\mi dz}-1}\bigg]\cosec(\pi z+\pi b_1)\cos(\pi z+\pi b_2).
\end{equation}

Using the same contour and the same procedure  as in Theorems 1 and 2, we get
\begin{eqnarray}\label{coscoseccossum}
&&\sum_{j=0}^{d-1} \cos\bigg(\frac{2\pi mj}{d}\bigg)\cosec\bigg(\frac{\pi j}{d}+\pi b_1\bigg)\cos\bigg(\frac{\pi j}{d}+\pi b_2\bigg)\nonumber\\&&=-d \cos[(2m-d)b_1\pi]\cosec (b_1d\pi)\cos[\pi+(b_2-b_1)\pi].
\end{eqnarray}

In the same way by replacing the last term by sine in (\ref{coscoseccosfun}), we get
\begin{eqnarray}\label{coscosecsinsum}
&&\sum_{j=0}^{d-1} \cos\bigg(\frac{2\pi mj}{d}\bigg)\cosec\bigg(\frac{\pi j}{d}+\pi b_1\bigg)\sin\bigg(\frac{\pi j}{d}+\pi b_2\bigg)\nonumber\\&&=-d \cos[(2m-d)b_1\pi]\cosec (b_1d\pi)\sin[\pi+(b_2-b_1)\pi].
\end{eqnarray}

Of course this sum, as before, is valid when $b_1\notin\mathbb{Z}/d$ such that the trigonometric sum is well defined.
By replacing the cosecant in (\ref{coscoseccosfun}) by secant, we get  following two sums

\begin{equation}\label{cosseccossum}
\sum_{j=0}^{d-1} \cos\bigg(\frac{2\pi mj}{d}\bigg)\sec\bigg(\frac{\pi j}{d}+\pi b_1\bigg)\cos\bigg(\frac{\pi j}{d}+\pi b_2\bigg)=
\begin{cases}(-1)^{m+1}d \cos[(2m-d)b_1\pi]\cosec (b_1d\pi)\\ \cos[\frac{\pi}{2}+(b_2-b_1)\pi] ; \qquad\qquad\qquad\qquad \text{if}\ d=\text{even}, \\
   (-1)^{m-d}d \sin[(2m-d)b_1\pi]\sec (b_1d\pi)\\ \cos[\frac{\pi}{2}+(b_2-b_1)\pi] ; \qquad\qquad\qquad\qquad  \text{if}\ d=\text{odd}.
    \end{cases}
\end{equation}

and

\begin{equation}\label{cossecsinsum}
\sum_{j=0}^{d-1} \cos\bigg(\frac{2\pi mj}{d}\bigg)\sec\bigg(\frac{\pi j}{d}+\pi b_1\bigg)\sin\bigg(\frac{\pi j}{d}+\pi b_2\bigg)=
\begin{cases}(-1)^{m+1}d \cos[(2m-d)b_1\pi]\cosec (b_1d\pi)\\ \sin[\frac{\pi}{2}+(b_2-b_1)\pi] ;\qquad\qquad\qquad\qquad  \text{if}\ d=\text{even}, \\
   (-1)^{m-d}d \sin[(2m-d)b_1\pi]\sec (b_1d\pi)\\ \sin[\frac{\pi}{2}+(b_2-b_1)\pi] ; \qquad\qquad\qquad\qquad \text{if}\ d=\text{odd}.
    \end{cases}
\end{equation}
For even $d$ the sums are well defined if $b_1\notin\mathbb{Z}/d$ and for odd $d$ the condition on $b_1$ is $b_1 \notin \mathbb{O}/2d$.

In the function  (\ref{coscoseccosfun}), if  we replace the minus sign between two exponentials by the plus sign, we will find
\begin{eqnarray}\label{sincoseccossum}
&&\sum_{j=1}^{d-1} \sin\bigg(\frac{2\pi mj}{d}\bigg)\cosec\bigg(\frac{\pi j}{d}+\pi b_1\bigg)\cos\bigg(\frac{\pi j}{d}+\pi b_2\bigg)\nonumber\\&&=d \sin[(2m-d)b_1\pi]\cosec (b_1d\pi)\cos[\pi+(b_2-b_1)\pi].
\end{eqnarray}

In the similar way replacing the last cosine term by sine, we get
\begin{eqnarray}\label{sincosecsinsum}
&&\sum_{j=1}^{d-1} \sin\bigg(\frac{2\pi mj}{d}\bigg)\cosec\bigg(\frac{\pi j}{d}+\pi b_1\bigg)\sin\bigg(\frac{\pi j}{d}+\pi b_2\bigg)\nonumber\\&&=d \sin[(2m-d)b_1\pi]\cosec (b_1d\pi)\sin[\pi+(b_2-b_1)\pi].
\end{eqnarray}
The condition on $b_1$ for the above two sums is $b_1\notin \mathbb{Z}/d$.

Just like we obtained (\ref{cosseccossum}) and (\ref{cossecsinsum}), we can find

\begin{equation}\label{sinseccossum}
\sum_{j=1}^{d-1} \sin\bigg(\frac{2\pi mj}{d}\bigg)\sec\bigg(\frac{\pi j}{d}+\pi b_1\bigg)\cos\bigg(\frac{\pi j}{d}+\pi b_2\bigg)=
\begin{cases} (-1)^{m}d \sin[(2m-d)b_1\pi]\cosec (b_1d\pi)\\ \cos[\frac{\pi}{2}+(b_2-b_1)\pi]; \qquad\qquad\qquad\qquad  \text{if}\ d=\text{even}, \\
  (-1)^{m-d}d \cos[(2m-d)b_1\pi]\sec (b_1d\pi)\\ \cos[\frac{\pi}{2}+(b_2-b_1)\pi] ; \qquad\qquad\qquad\qquad  \text{if}\ d=\text{odd}.
    \end{cases}
\end{equation}

and

\begin{equation}\label{sinsecsinsum}
\sum_{j=1}^{d-1} \sin\bigg(\frac{2\pi mj}{d}\bigg)\sec\bigg(\frac{\pi j}{d}+\pi b_1\bigg)\sin\bigg(\frac{\pi j}{d}+\pi b_2\bigg)=
\begin{cases} (-1)^{m}d \sin[(2m-d)b_1\pi]\cosec (b_1d\pi)\\ \sin[\frac{\pi}{2}+(b_2-b_1)\pi]; \qquad\qquad\qquad\qquad  \text{if}\ d=\text{even}, \\
  (-1)^{m-d}d \cos[(2m-d)b_1\pi]\sec (b_1d\pi)\\ \sin[\frac{\pi}{2}+(b_2-b_1)\pi] ; \qquad\qquad\qquad\qquad  \text{if}\ d=\text{odd}.
    \end{cases}
\end{equation}
For even $d$ the above two sums are well defined if $b_1\notin\mathbb{Z}/d$ and for odd $d$ the condition is $b_1\notin\mathbb{O}/2d$.

Let us now take a different complex function. We can then obtain many more sums.
\begin{equation}\label{coscoseccosecfun}
f(z)=\bigg[\frac{\me^{2\pi\mi mz }}{\me^{2\pi\mi dz}-1}-\frac{\me^{-2\pi\mi mz }}{\me^{-2\pi\mi dz}-1}\bigg]\cosec(\pi z+\pi b_1)\cosec(\pi z+\pi b_2).
\end{equation}

By using residue theorem, we can show that
\begin{eqnarray}\label{coscoseccosecsum}
&&\sum_{j=0}^{d-1} \cos\bigg(\frac{2\pi mj}{d}\bigg)\cosec\bigg(\frac{\pi j}{d}+\pi b_1\bigg)\cosec\bigg(\frac{\pi j}{d}+\pi b_2\bigg)\\&&=-d \cos[(2m-d)b_1\pi]\cosec (b_1d\pi)\cosec[\pi+(b_2-b_1)\pi] \nonumber\\&&\quad-d \cos[(2m-d)b_2\pi]\cosec (b_2d\pi)\cosec[\pi+(b_1-b_2)\pi].\nonumber
\end{eqnarray}
 This sum is valid for $b_1\notin\mathbb{Z}/d$ and $b_2\notin\mathbb{Z}/d$ such that the trigonometric sum is well defined.
In (\ref{coscoseccosecfun}) instead of two cosecant, if we use one cosecant and one secant, we will get

\begin{equation}\label{coscosecsecsum}
\sum_{j=0}^{d-1} \cos\bigg(\frac{2\pi mj}{d}\bigg)\cosec\bigg(\frac{\pi j}{d}+\pi b_1\bigg)\sec\bigg(\frac{\pi j}{d}+\pi b_2\bigg)=
\begin{cases} \Big(-d \cos[(2m-d)b_1\pi]\cosec (b_1d\pi)\\\sec[\pi+(b_2-b_1)\pi]+(-1)^{m+1}d \cos[(2m-d)b_2\pi]\\\cosec (b_2d\pi) \cosec[\frac{\pi}{2}+(b_1-b_2)\pi]\Big); \quad \text{if}\ d=\text{even}, \\
  \Big(-d \cos[(2m-d)b_1\pi]\cosec (b_1d\pi)\\\sec[\pi+(b_2-b_1)\pi] +(-1)^{m-d}d \sin[(2m-d)b_2\pi]\\\sec (b_2d\pi) \cosec[\frac{\pi}{2}+(b_1-b_2)\pi]\Big) ; \qquad \text{if}\ d=\text{odd}.
    \end{cases}
\end{equation}
For even $d$ the above sum is valid for $b_1\notin\mathbb{Z}/d$ and $b_2\notin\mathbb{Z}/d$. In case of odd $d$ the conditions are $b_1\notin\mathbb{Z}/d$ and $b_2\notin\mathbb{O}/2d$.
Again if both are secants, then

\begin{equation}\label{cossecsecsum}
\sum_{j=0}^{d-1} \cos\bigg(\frac{2\pi mj}{d}\bigg)\sec\bigg(\frac{\pi j}{d}+\pi b_1\bigg)\sec\bigg(\frac{\pi j}{d}+\pi b_2\bigg)=
\begin{cases} \Big((-1)^{m+1}d \cos[(2m-d)b_1\pi]\cosec (b_1d\pi)\\ \sec[\frac{\pi}{2}+(b_2-b_1)\pi] +(-1)^{m+1}d \cos[(2m-d)b_2\pi]\\ \cosec (b_2d\pi) \sec[\frac{\pi}{2}+(b_1-b_2)\pi]\Big); \qquad \text{if}\ d=\text{even}, \\
  \Big((-1)^{m-d}d \sin[(2m-d)b_1\pi]\sec (b_1d\pi)\\ \sec[\frac{\pi}{2}+(b_2-b_1)\pi]+(-1)^{m-d}d \sin[(2m-d)b_2\pi]\\ \sec (b_2d\pi)\sec[\frac{\pi}{2}+(b_1-b_2)\pi]\Big) ; \qquad\quad \text{if}\ d=\text{odd}.
    \end{cases}
\end{equation}

The conditions on $b$'s for even $d$ are $b_1\notin\mathbb{Z}/d$ and $b_2\notin\mathbb{Z}/d$. For odd $d$ case the conditions are 
$b_1\notin\mathbb{O}/2d$ and $b_2\notin\mathbb{O}/2d$.

Like (\ref{coscoseccosecsum}), (\ref{coscosecsecsum}) and (\ref{cossecsecsum}), we can get similar sums with sine as follows.

\begin{eqnarray}\label{sincoseccosecsum}
&&\sum_{j=1}^{d-1} \sin\bigg(\frac{2\pi mj}{d}\bigg)\cosec\bigg(\frac{\pi j}{d}+\pi b_1\bigg)\cosec\bigg(\frac{\pi j}{d}+\pi b_2\bigg)\\&&=d \sin[(2m-d)b_1\pi]\cosec (b_1d\pi)\cosec[\pi+(b_2-b_1)\pi] \nonumber\\&&\quad+d \sin[(2m-d)b_2\pi]\cosec (b_2d\pi)\cosec[\pi+(b_1-b_2)\pi].\nonumber
\end{eqnarray}
Conditions on $b_1$ and $b_2$ are $b_1\notin\mathbb{Z}/d$ and $b_2\notin\mathbb{Z}/d$.

\begin{equation}\label{sincosecsecsum}
\sum_{j=1}^{d-1} \sin\bigg(\frac{2\pi mj}{d}\bigg)\cosec\bigg(\frac{\pi j}{d}+\pi b_1\bigg)\sec\bigg(\frac{\pi j}{d}+\pi b_2\bigg)=
\begin{cases} \Big(d \sin[(2m-d)b_1\pi]\cosec (b_1d\pi)\\\sec[\pi+(b_2-b_1)\pi]+(-1)^{m}d \sin[(2m-d)b_2\pi]\\\cosec (b_2d\pi) \cosec[\frac{\pi}{2}+(b_1-b_2)\pi]\Big); \quad \text{if}\ d=\text{even}, \\
  \Big(d \sin[(2m-d)b_1\pi]\cosec (b_1d\pi)\\\sec[\pi+(b_2-b_1)\pi] +(-1)^{m-d}d \cos[(2m-d)b_2\pi]\\\sec (b_2d\pi) \cosec[\frac{\pi}{2}+(b_1-b_2)\pi]\Big) ; \qquad \text{if}\ d=\text{odd}.
    \end{cases}
\end{equation}
For even $d$ the above sum is valid for $b_1\notin\mathbb{Z}/d$ and $b_2\notin\mathbb{Z}/d$. In case of odd $d$ the conditions are $b_1\notin\mathbb{Z}/d$ and $b_2\notin\mathbb{O}/2d$.


\begin{equation}\label{sinsecsecsum}
\sum_{j=0}^{d-1} \sin\bigg(\frac{2\pi mj}{d}\bigg)\sec\bigg(\frac{\pi j}{d}+\pi b_1\bigg)\sec\bigg(\frac{\pi j}{d}+\pi b_2\bigg)=
\begin{cases} \Big((-1)^{m}d \sin[(2m-d)b_1\pi]\cosec (b_1d\pi)\\\sec[\frac{\pi}{2}+(b_2-b_1)\pi] +(-1)^{m}d \sin[(2m-d)b_2\pi]\\\cosec (b_2d\pi) \sec[\frac{\pi}{2}+(b_1-b_2)\pi]\Big); \qquad \text{if}\ d=\text{even}, \\
  \Big((-1)^{m-d}d \cos[(2m-d)b_1\pi]\sec (b_1d\pi)\\\sec[\frac{\pi}{2}+(b_2-b_1)\pi]+(-1)^{m-d}d \cos[(2m-d)b_2\pi] \\\sec (b_2d\pi) \sec[\frac{\pi}{2}+(b_1-b_2)\pi]\Big) ; \quad\qquad \text{if}\ d=\text{odd}.
    \end{cases}
\end{equation}
The conditions on $b$'s for even $d$ are $b_1\notin\mathbb{Z}/d$ and $b_2\notin\mathbb{Z}/d$. For odd $d$ case the conditions are 
$b_1\notin\mathbb{O}/2d$ and $b_2\notin\mathbb{O}/2d$.
  
   In this section, we have considered the product of three trigonometric functions with simple powers. One can extend these results to higher powers
   of trigonometric functions, as well as products of more than three trigonometric functions. The list is endless. We have only illustrated a few cases.

\section{Conclusion}
   
     We have obtained a number of finite sums involving products of two
   or more trigonometric functions. They were mostly based on a specific
   choice of contour and a wide variety of integrands. Many more sums can
   be obtained. For the simpler powers of trigonometric functions, these
   sums can be remarkably simple. Most of these sums involve tangents, 
   cotangents, secants, and cosecants. For these functions, we don't have
   simple expansions for the sum and difference of variables in their arguments. Therefore,
   one has to compute them independently. We have calculated the sums 
   for a few cases. One can extend these calculations in many different directions, as discussed
   in the above sections. Most of the simpler expressions
   should be in handbooks that have trigonometric sums.

\end{document}